\documentclass[10pt,reqno]{amsart}
\usepackage{epsfig}
\usepackage{color}

\usepackage{amsmath}
\usepackage{esint}
\usepackage[colorlinks=true]{hyperref}
\usepackage{paralist}

\usepackage{mathrsfs}

\allowdisplaybreaks
\numberwithin{equation}{section}

\def\H{{\mathcal H}}
\def\L{{\mathcal L}}
\def\N{{\bf N}}
\def\R{{\bf R}}

\def\Z{{\bf Z}}

\def\={:=}
\def\eps{{\varepsilon}_j}
\def\epsN{\eps^{N-1}}
\def\HH{{\H}^{N-1}}

\def\FF{{\mathcal F}}
\def\FFepsi{{\mathcal F}_{j}}

\def\ii{{\tt i}}
\def\kk{{\tt k}}
\def\e{\varepsilon}

\def\Wsp{W^{1,p}}
\def\KpmuU{K^p(U,\mu)}
\def\Kpmu{K^p(\R^N,\mu)}

\def\Wspo{W_0^{1,p}}
\def\WspoG{W_{0,\Sigma}^{1,p}}

\def\capmu{\mathrm{cap}_{p,\mu}}
\def\Om{\Omega}
\def\om{\omega}
\def\deltaj{\delta_j}
\def\deltae{\delta_{j}}

\def\Qij{Q_j^\ii}
\def\Qijo{Q_j^\ii}
\def\Cij{C_j^{\ii,h_\ii}}

\def\Bji{B_j^{\ii}}

\def\Bijh{B_j^{\ii,h}}

\def\Cijh{C_j^{\ii,h}}
\def\Sijh{S_j^{\ii,h}}

\def\res{\mathop{\hbox{\vrule height 7pt width .4pt depth 0pt
\vrule height .4pt width 6pt depth 0pt}}\nolimits}

\def\Teioe{T_{j}^\ii(\om)}
\def\Teio{T_{j}^\ii(\om)}
\def\Teiot{\tilde{T}_{j}^\ii(\om)}
\def\Teoe{T_{j}(\om)}
\def\Teo{T_{j}(\om)}
\def\gammaio{\gamma(\ii,\om)}
\def\Ieps{{\mathcal I}_j}
\def\IIeps{{\mathscr I}_j}
\def\gammab{\gamma_0}
\def\xijo{y_j^{\ii}}
\def\xije{x_{j}^\ii}
\def\xiij{\xi_j^\ii}
\def\tildeu{\tilde{u}}
\def\FFjpsi{{\mathcal F}_j}
\def\FFpsi{{\mathcal F}}
\def\vj{v_j}
\def\mmu{\mathbb{P}}
\def\xn{x_N}
\def\xnc{y}

\def\UU{\partial_NU}
\def\LN{\L^{N-1}}
\def\lambdaj{\lambda_j}
\def\Ur{U_r}

\def\res{\mathop{\hbox{\vrule height 7pt width .4pt depth 0pt
\vrule height .4pt width 6pt depth 0pt}}\nolimits}

\def \trait (#1) (#2) (#3){\vrule width #1pt height #2pt depth #3pt}
\def \qed{\hfill
        \trait (0.1) (6) (0)
        \trait (6) (0.1) (0)
        \kern-6pt
        \trait (6) (6) (-5.9)
        \trait (0.1) (6) (0)
\medskip}

\newtheorem{theorem}{Theorem}[section]
\newtheorem{lemma}[theorem]{Lemma}
\newtheorem{corollary}[theorem]{Corollary}
\newtheorem{remark}[theorem]{Remark}

\newtheorem{proposition}[theorem]{Proposition}

\title[{$\Gamma$-convergence for random fractional obstacles}]
{Homogenization of random fractional obstacle problems via
$\Gamma$-convergence}

\begin{document}
\maketitle
\smallskip
\begin{center}
\begin{tabular}{c@{\qquad}c@{}}
Matteo Focardi \\[-0.3cm]
{\footnotesize Dip. Mat. ``U. Dini'' }
\\[-0.3cm]
{\footnotesize V.le Morgagni, 67/a}
\\[-0.3cm]
{\footnotesize I-50134  Firenze}
{\footnotesize focardi@math.unifi.it}
\end{tabular}
\end{center}
\medskip
\begin{abstract}
$\Gamma$-convergence methods are used to prove homogenization
results for fractional obstacle problems in periodically perforated 
domains.
The obstacles have random sizes and shapes and their
capacity scales according to a stationary ergodic process.
We use a trace-like representation of fractional
Sobolev norms in terms of weighted Sobolev energies established
in \cite{Caf-Silv}, a weighted ergodic theorem and a joining
lemma in varying domains following the approach by \cite{ANB}.

Our proof is alternative to those contained in \cite{Caf-Mel1},
\cite{Caf-Mel2}.

\end{abstract}


\section{Introduction}\label{intro}

The homogenization of (non-)linear elliptic obstacle problems in
periodically perforated domains has received much attention after
the seminal papers of Marchenko and Khruslov~\cite{MK}, Rauch and
Taylor~\cite{RT1},\cite{RT2} and Cioranescu and Murat~\cite{CM} (see
\cite{CC},\cite{AP},\cite{DGDML},\cite{DML},\cite{DM1},\cite{DM2},\cite{DM3},\cite{LOPS},\cite{ANB} and \cite{ATT},\cite{BDF},\cite{CD},\cite{DM}
for a more exhaustive list of references).
The problem has been successfully tackled by making use of abstract 
techniques of $\Gamma$-convergence, and fully solved in a series of papers
by Dal Maso \cite{DM1},\cite{DM2},\cite{DM3}. 
A constructive approach in the periodic case for bilateral obstacles has
been developed by Ansini and Braides \cite{ANB}.
In general, a relaxation process takes place and the limit problem
contains a finite penalization term related to the capacity
of the homogenizing obstacles.

All the quoted results deal with Sobolev type energies and deterministic
distributions of the set of obstacles with deterministic sizes and shapes.
More recently, two papers \cite{Caf-Mel1}, \cite{Caf-Mel2} have enlarged
the stage to fractional Sobolev energies and by considering random
sizes and shapes for the obstacles. 
More precisely, given a probability
space $(\Om,\mathscr{P},\mmu)$, for all $\om\in\Om$ consider a periodic
distribution of sets $T_\e(\om)$ and let $v_\e(\cdot,\om)$ be the solution
of the problem
\begin{eqnarray}
  \label{eq:pbm}
  \begin{cases}
(-\triangle)^sv(y)\geq 0 & y\in\R^{N-1}\\
(-\triangle)^sv(y)= 0 & y\in\R^{N-1}\setminus T_\e(\om), \text{ and }
y\in T_\e(\om) \text{ if } v(y)>\psi(y)\\
v(y)\geq\psi(y) & y\in T_\e(\om).
\end{cases}
\end{eqnarray}
The operator $(-\triangle)^s$ is the fractional Laplace operator of
order $s\in(0,1)$ defined in terms of the Fourier transform, by
$\mathscr{F}((-\triangle)^sv)(\xi)=|\xi|^{2s}\hat{v}(\xi)$; $\psi$
is the \emph{obstacle function} and it is assumed to be in $C^{1,1}(\R^{N-1})$.
In case $s=1/2$ the minimum problem in \eqref{eq:pbm} is known as 
Signorini's problem and it is related to a semi-permeable membrane model. 
We refer to the papers \cite{Caf-Mel1} and \cite{Caf-Mel2} for a 
more detailed description of the underlying physical model.

Problem \eqref{eq:pbm} has a natural variational character.
Indeed, it can be interpreted as the Euler-Lagrange
equation solved by the minimizer of 
\begin{eqnarray}\label{eq:pbm1}
\inf_{\dot{H}^s(\R^{N-1})}\left\{\|v\|_{\dot{H}^s(\R^{N-1})}:\,v\geq\psi
\text{ on } T_\e(\om)\right\}.
\end{eqnarray}
Here $\|v\|_{\dot{H}^s(\R^{N-1})}=
\||\xi|^{2s}\hat{v}(\xi)\|_{L^2(\R^{N-1})}$ is the usual
norm in the homogeneous fractional Sobolev space $\dot{H}^s(\R^{N-1})$.

An additional variational characterization of problem \eqref{eq:pbm}
can be given by following the work by Caffarelli and
Silvestre~\cite{Caf-Silv} who have represented fractional Sobolev
norms on $\R^{N-1}$ in terms of boundary value problems for
degenerate (but local!) elliptic equations in the higher dimensional
half-space $\R^N_+$; equivalently, in terms of minimal energy extensions
of a (suitable) weighted Dirichlet integral as for the harmonic extension
of $\dot{H}^{1/2}$ functions.
It turns then out that the extension $u_\e(\cdot,\om)$ of
$v_\e(\cdot,\om)$ to $\R^N_+$ solves the problem
\begin{eqnarray}
  \label{eq:pbm2}
\inf_{W^{1,2}(\R^N_+,|\xn|^a)}\left\{\int_{\R^N_+}|\xn|^a|\nabla u(y,\xn)|^2d\L^N:\,
u(y,0)\geq\psi(y)\,\,y\in T_\e(\om)\right\}.
\end{eqnarray}
Here, the parameter $a$ ruling the degeneracy of the elliptic
equation equals $2s-1$ (and thus belongs to $(-1,1)$), and
$W^{1,2}(\R^N_+,|\xn|^a)$ is the weighted Sobolev space associated
to the measure $|\xn|^ad\L^N(y,\xn)$.

To investigate the asymptotic behaviour of $u_\e(\cdot,\om)$ as
$\e$ vanishes some assumptions have to be imposed on the obstacles set
$T_\e(\om)$. Mild hypotheses have been introduced in
\cite{Caf-Mel1}, \cite{Caf-Mel2}: the set $T_\e(\om)$ is the union
of periodically distributed sets (but with random sizes and shapes!)
whose capacity scales according to a stationary and ergodic process
$\gamma$ (see (Hp 1) and (Hp 2) in Section~\ref{assmptns}).
Under these assumptions Caffarelli and Mellet~\cite{Caf-Mel2}
have proven that there exists a constant $\alpha_0\geq 0$ such that
the solution $u_\e(\cdot,\om)$ of \eqref{eq:pbm2} converges locally weakly in
$W^{1,2}(\R^N_+,|\xn|^a)$ and $\mmu$ a.s. in $\Om$ to the solution
$\bar{u}$ of
\begin{eqnarray}
  \label{eq:pbm3}
  \inf_{W^{1,2}(\R^N_+,|\xn|^a)}\left\{\int_{\R^N_+}|\xn|^a|\nabla u(y,\xn)|^2d\L^N+
\frac {\alpha_0}2\int_{\R^{N-1}}|(\psi(y)-u(y,0))\vee 0|^2dy\right\}.
\end{eqnarray}
The proof of such a result relies
on the regularity of fractional obstacle problems established
by Caffarelli et al. \cite{Caf-Sal-Silv}, and on the PDEs approach
to homogenization  based on the Tartar's oscillating test function method
(see \cite{T}, \cite{CM} and \cite{CD} for further references).

The aim of this paper is to give an alternative elementary proof
of the above quoted homogenization results via $\Gamma$-convergence
techniques.
We are able to avoid the use of the regularity
theory developed in \cite{Caf-Sal-Silv} and thus to relax the
smoothness assumption on the obstacle function $\psi$.
In addition, we determine explicitely the constant $\alpha_0$ in the
capacitary contribution of the limit energy, and show that it equals
the expectation of the process $\mathbb{E}[\gamma]$ (see Theorem~\ref{main}).

Despite this, the proof is not self-contained since we still use
the trace-like representation for fractional norms established
in \cite{Caf-Silv}.
A direct approach is still under investigation, and deserves
additional efforts since the difficulties introduced in the problem
by the non-locality of fractional energies.

The main tools of our analysis are a joining lemma in varying boundary
domains for weighted energies and a weighted version of Birkhoff's
ergodic theorem. The joining lemma follows the line of
the analogous result in perforated open sets
for standard Sobolev spaces proved by Ansini and Braides~\cite{ANB}.
It is a variant of an idea by De Giorgi \cite{DeG} 
in the setting of varying domains, on the way of
matching boundary conditions by increasing the energy only 
up to a small error. This method is elementary and based
on a clever slicing and averaging argument, looking for those zones where
the energy does not concentrate.
The joining lemma allows us to reduce in the $\Gamma$-limit process to
families of functions which are constants on suitable annuli surrounding
the obstacle sets.
Thus, to estimate the capacitary contribution close to the
obstacle set $T_\e(\om)$ we exploit the capacitary scaling
assumption on the process $\gamma$ together with
a weighted variant of Birkhoff's ergodic theorem (see Theorem~\ref{ergo-w}).
This argument allows us to show that $\alpha_0$ equals $\mathbb{E}[\gamma]$.

In Section~\ref{assmptns} we list the assumptions
and state the homogenization result. To avoid unnecessary generality
we deal with the model case of $p$-norms, $p\in(1,+\infty)$, since 
this case contains all the features of the problem.
Section~\ref{sobolev} collects several results concerning weighted
Sobolev spaces in case the weight function is a Muckenhoupt weight
of the form $w(y,\xn)=|\xn|^a$.
A weighted ergodic theorem relevant in our analysis is
proved in Section~\ref{weighted-ergodic}.
In Section~\ref{mainresult} we prove the $\Gamma$-convergence theorem.
Finally, in Section~\ref{genztn} we indicate several possible
generalizations. 

\section{Statement of the Main Result}\label{assmptns}

\subsection{Basic Notations}
The ball in $\R^N$ with centre $x$ and radius $r>0$ is denoted
by $B_r(x)$, and simply by $B_r$ in case $x=\underline{0}$.
The interior and the closure of a set $E\subset\R^N$ are denoted 
by $\mathrm{int(E)}$ and $\overline{E}$, respectively.
Given two sets $E\subset\subset F$ in $\R^N$, a 
\emph{cut-off function between E and F} is any 
$\varphi\in C^\infty_0(F)$ such that $\varphi|_E\equiv 1$. 

Not to overburden the notation each set $E\subseteq \R^{N-1}$ and
its copy $E\times\{0\}\subseteq\R^N$ will be undistinguished.

In the sequel $U$ denotes any connected open subset of the half-space
$\R^N_+:=\{x=(\xnc,\xn):\,\xnc\in\R^{N-1},\, \xn>0\}$ whose boundary is
Lipschitz regular. The part of the boundary of $U\subseteq\R^N_+$
lying on $\{\xn=0\}$ is denoted by $\UU:=\partial U\cap\{\xn=0\}$.

We use standard notations for Hausdorff and Lebesgue measures,
and Lebesgue spaces.
The integration with respect to the measure $\H^{N-1}\res\{\xn=0\}$ is
denoted by $d\xnc$, and for $V\subseteq\{\xn=0\}$
the spaces $L^p(V,\H^{N-1}\res\{\xn=0\})$
simply by $L^p(V)$, $p\in[1,+\infty]$.

The lattice in $\R^{N-1}$ underlying the periodic homogenization
process is identified via the points $\xijo:=\ii\eps\in\R^{N-1}$, 
$\xije:=(\xijo,0)\in\R^N$ and the cubes
$\Qijo:=\xijo+\eps[-1/2,1/2)^{N-1}\subset\R^{N-1}$, $\ii\in\Z^{N-1}$.
Here, $(\eps)_j$ is a positive infinitesimal sequence.
Finally, for any set $E\subseteq \R^{N-1}$ define
$$
\Ieps(E):=\{\ii\in\Z^{N-1}:\,\Qijo\subseteq E\}. 
$$

\subsection{$\Gamma$-convergence}
We recall the notion of $\Gamma$-convergence introduced by De Giorgi
in a generic metric space $(X,d)$ endowed with the topology induced
by $d$ (see \cite{DM},\cite{B3}).
A sequence of functionals $F_j:X\to [0,+\infty]$
{\it $\Gamma$-converges} to a functional $F:X\to [0,+\infty]$ in $u\in X$,
in short $F(u)=\Gamma\hbox{-}\lim_{j}F_j(u)$,
if 
the following two conditions hold:

{(i)} ({\it liminf inequality}) $\forall\  (u_j)$
converging to $u$ in $X$,
we have $\liminf_jF_{j}(u_j)\ge F(u)$;

{(ii)} ({\it limsup inequality}) $\exists$ $(u_j)$ converging to $u$ in $X$
such that $\limsup_jF_{j}(u_j)\le F(u)$.

\noindent We say that $F_j$ \emph{$\Gamma$-converges} to $F$
(or $F$= $\Gamma$-lim$_{j}F_j$) if
$F(u)=\Gamma\hbox{-}\lim_{j}F_j(u)$ $\forall u\in X$.
We may also define the \emph{lower} and \emph{upper $\Gamma$-limits} as
$$
\Gamma\hbox{-}\limsup_{j}F_j(u)=\inf\{\limsup_{j}
F_j(u_j):\  u_j\to u\},
$$
$$
\Gamma\hbox{-}\liminf_{j}F_j(u)=\inf\{\liminf_{j}
F_j(u_j):\  u_j\to u\},
$$
respectively, so that conditions (i) and (ii) are equivalent to
$\Gamma$-limsup$_j F_j(u)=\Gamma$-liminf$_j F_j(u)=F(u)$.
Moreover, the functions $\Gamma$-limsup$_j F_j$ and
$\Gamma$-liminf$_j F_j$ are lower semicontinuous.

One of the main reasons for the introduction of this notion is explained
by the following fundamental theorem (see \cite[Theorem 7.8]{DM}).

\begin{theorem}\label{min}
Let $F=\Gamma$-$\lim_{j}F_j$, and assume there exists
a compact set $K\subset X$  such that
$\inf_X F_j=\inf_K F_j$ for all $j$. Then there exists
$\min_X F =\lim_{j}\inf_X F_j$. Moreover,
if $(u_j)$ is a converging sequence such that
$\lim_j F_{j}(u_j)=\lim_j\inf_X F_{j}$
then its limit is a minimum point for $F$.
\end{theorem}

\subsection{Assumptions and Statement of the Main Result}

We consider a probability space $(\Om,\mathscr{P},\mmu)$.
For all $\om\in\Om$ and $j\in\N$ the set $\Teo\subseteq\R^{N-1}$
is given by
$$
\Teoe=\cup_{\ii\in\Z^{N-1}}\Teioe
$$
where the sets $\Teioe\subseteq\Qij$ satisfy the following conditions:
\begin{itemize}
\item[{\bf (Hp 1).}] \emph{Capacitary Scaling:}
There exist a positive infinitesimal sequence $(\deltaj)_j$ and
a process $\gamma:\Z^{N-1}\times\Om\to[0,+\infty)$ such that
for all $\ii\in\Z^{N-1}$ and $\om\in\Om$
$$
\capmu(\Teioe)=\deltae\gammaio.
$$

\item[{\bf (Hp 2).}] \emph{Ergodicity $\&$ Stationarity of the Process:}
The process $\gamma:\Z^{N-1}\times\Om\to[0,+\infty)$ is stationary ergodic:
There exists a family of measure-preserving transformations
$\tau_\kk:\Om\to\Om$ satisfying for all $\ii,\kk\in\Z^{N-1}$ and $\om\in\Om$
\begin{equation}
  \label{eq:ergogamma}
\gamma(\ii+\kk,\om)=\gamma(\ii,\tau_\kk\om),
\end{equation}
and such that if $A\subseteq\Om$ is an invariant set, i.e. $\tau_\kk A=A$
for all $\kk\in\Z^{N-1}$, then either $\mmu(A)=0$ or $\mmu(A)=1$.

Moreover, for some $\gammab>0$ we have for all $\ii\in\Z^{N-1}$ and
$\mmu$ a.s. $\om\in\Om$
$$
\gammaio\leq\gammab.
$$

\item[{\bf (Hp 3).}] \emph{Strong Separation:}
There exist $\e$, $M>0$ such that for all $\ii\in\Z^{N-1}$, $\om\in\Om$,
and for every $\eps\in(0,\e)$ it holds
$\Teioe\subseteq \xijo+M\eps^{\beta}[-1/2,1/2)^{N-1}$, where
$\beta=(N-1)/(N-p+a)$.

\item[{\bf (Hp 4).}] The sequence $(\deltae\eps^{-N+1})$ has a limit in
$[0,+\infty]$. We denote such a value $\Lambda$.

\end{itemize}
Assumptions (Hp 1)-(Hp 3) were introduced in \cite{Caf-Mel2} 
(see Remark~\ref{hptre} for a weak variant of (Hp 3)).

In the following remarks we briefly comment on the previous assumptions.

\begin{remark}
The capacitary scaling assumption implies that
$$
\capmu\left(\Teoe\right)\leq\sum_{\ii\in\Z^{N-1}}\capmu(\Teioe)=
\deltae\sum_{\ii\in\Z^{N-1}}\gammaio.
$$
Heuristically, 
we may assume $\capmu\left(\Teoe\right)\sim\sum_{\ii\in\Z^{N-1}}\capmu(\Teioe)$
since the obstacles $\Teioe$ are sufficiently far apart one from
the other by the strong separation assumption.
Hence, by taking into account Birkhoff's individual ergodic theorem
$\mmu$ a.s. in $\Om$ we infer
$$
\capmu\left(\Teoe\right)\sim\Lambda\mathbb{E}[\gamma].
$$
Thus we can distinguish three regimes according to the asymptotic
behaviour of $\deltae\e^{-N+1}$ (see Theorem~\ref{main}).
\end{remark}

\begin{remark}
The stationarity 
property is a mild assumption in order to have some averaging along the
homogenization process, a condition weaker than periodicity or 
quasi-periodicity. 
It implies that the random field $\gamma$ is statistically 
homogeneous w.r.to the action of traslations compatible with the 
underlying periodic lattice, e.g. the random variables 
$\gamma(\ii,\cdot)$ are independent and identically distributed. 
\end{remark}

With fixed exponents $a\in(-1,+\infty)$ and $p\in((1+a)\vee 1,N+a)$
(these restrictions will be justified in Section~\ref{sobolev},
Remark~\ref{trivial} and Appendix~\ref{Muck}), consider the measure
$\mu:=|\xn|^ad\L^N$ and the corresponding weighted Sobolev space
$\Wsp(U,\mu)$ (see Section~\ref{sobolev}).

Let $\psi$ be upper bounded and continuous in the relative interior of $\UU$
w.r.to the relative topology of $\{\xn=0\}$ (for some comments on this 
assumption see Remark~\ref{uncont}) and define the functional 
$\FFepsi:L^p(U,\mu)\times\Om\to[0,+\infty]$ by
\begin{eqnarray}  \label{eq:fapprox}
  \FFepsi(u,\om)=
  \begin{cases}
\displaystyle{\int_U|\nabla u|^p\,d\mu} &
\text{ if } u\in \Wsp(U,\mu),\,
\tildeu\geq\psi\,\,  \capmu \text{ q.e. on } \Teoe\cap\UU\\
+\infty & \text{ otherwise. }
  \end{cases}
\end{eqnarray}
Here, $\capmu$ is the variational $(p,\mu)$-capacity associated with $\mu$,
and $\tilde{u}$ denotes the precise representative of $u$ which
is defined except on a set of capacity zero (see Section~\ref{sobolev}).

To state the main result of the paper and not to make it trivial
we also assume that (see Remark~\ref{uncont1})
\begin{itemize}
\item[{\bf (Hp 5).}] there exists $f\in\Wsp(U,\mu)$ such that
$\tilde{f}\geq\psi$ $\capmu$ q.e. on $\UU$.
\end{itemize}

\begin{theorem}\label{main}
Assume (Hp 1)-(Hp 5) hold true, $N\geq 2$, and that $a\in(-1,+\infty)$,
$p\in((1+a)\vee 1,N+a)$.

Then there exists a set $\Om^\prime\subseteq\Om$ of full probability
such that for all $\om\in\Om^\prime$ the sequence
$(\FFepsi(\cdot,\om))_j$ $\Gamma$-converges in the $L^p(U,\mu)$
topology to the functional $\FFpsi:L^p(U,\mu)\to[0,+\infty]$ defined by
\begin{equation}\label{eq:Glimit}
\FFpsi(u)=\int_U|\nabla u|^p\,d\mu+\frac 12
\Lambda\mathbb{E}[\gamma]
\int_{\UU}|(\psi(\xnc)-u(\xnc,0))\vee 0|^p\,d\xnc
\end{equation}
if $u\in \Wsp(U,\mu)$, $+\infty$ otherwise.
\end{theorem}
In case $U$ is not bounded equi-coercivity for the functionals 
$\FFepsi$ is ensured only in the $L^p_{\mathrm{loc}}(U,\mu)$ topology. 
A relaxation phenomenon takes place and the domain of the limit 
has to be slightly enlarged according to Sobolev-Gagliardo-Nirenberg 
inequality in 
\begin{equation}\label{kpmu}
\KpmuU=\{u\in L^{p^\ast}(U,\mu): \nabla u\in (L^p(U,\mu))^N\},
\end{equation}
where $p^\ast=(N+a)p/(N+a-p)$ is the Sobolev exponent relative 
to $\Wsp(\R^N,\mu)$ (see Lemma~\ref{Gagliardo}). 
We show $\Gamma$-convergence in that case, too
\begin{theorem}\label{main-unbdd}
Under the assumptions of Theorem~\ref{main}, if $U$ is unbounded 
there exists a set $\Om^\prime\subseteq\Om$ of full probability
such that for all $\om\in\Om^\prime$
the family $(\FFepsi(\cdot,\om))_j$ $\Gamma$-converges in the
$L^p_{\mathrm{loc}}(U,\mu)$ topology to  the functional 
 $\FFpsi:L^p_{\mathrm{loc}}(U,\mu)\to[0,+\infty]$ defined by
\begin{equation}\label{eq:Glimit}
\FFpsi(u)=\int_U|\nabla u|^p\,d\mu+\frac 12
\Lambda\mathbb{E}[\gamma]
\int_{\UU}|(\psi(\xnc)-u(\xnc,0))\vee 0|^p\,d\xnc
\end{equation}
if $u\in\KpmuU$, $+\infty$ otherwise.

\end{theorem}
The set $\Om^\prime$ referred to in the statements of Theorem~\ref{main},
\ref{main-unbdd} is defined in Section~\ref{mainresult} below.

Theorem~\ref{main} is compatible with the addition of boundary data.
Assume that $U$ is bounded, denote by $\Sigma$ a non-empty and 
relatively open subset of $\partial U\setminus\UU$, and by 
$\WspoG(U,\mu)$ the strong closure in $\Wsp(U,\mu)$ of the restrictions 
to $U$ of functions $C^\infty(\R^N)$ vanishing on a neighbourhood of 
$\overline{\Sigma}$. Further, we require that 
$\overline{\Sigma}\cap\UU=\emptyset$ to avoid additional technicalities.
\begin{corollary}\label{bdryvalue}
Assume that $U$ is bounded, and that (Hp 1)-(Hp 4) hold true. 
With fixed $N\geq 2$, $a\in(-1,+\infty)$, $p\in((1+a)\vee 1,N+a)$ and 
$u_0\in\Wsp(U,\mu)$ s.t. $\tilde{u}_0\geq\psi$ $\capmu$ q.e.
on $\UU$ there exists a set $\Om^{\prime\prime}\subseteq\Om$ 
of full probability such that for all $\om\in\Om^{\prime\prime}$
the functionals $\FFepsi(\cdot,\om)+\mathscr{X}_{u_0+\WspoG(U,\mu)}$
$\Gamma$-converge in the $L^p(U,\mu)$ topology
to $\FFpsi+\mathscr{X}_{u_0+\WspoG(U,\mu)}$,
where $\mathscr{X}_{u_0+\WspoG(U,\mu)}$ is the $0,+\infty$
characteristic funtion of the subspace $u_0+\WspoG(U,\mu)$.
\end{corollary}
$\Gamma$-convergence theory then implies convergence of minimizers
provided the equi-coercivity of the $\FFepsi$'s holds 
(see Theorem~\ref{min}).
That property is ensured by Theorem~8 \cite{K2} in case $U$ is 
bounded, and by Lemma~\ref{Gagliardo} 
below if $U$ is unbounded.

\begin{corollary}
Under the assumptions of Corollary~\ref{bdryvalue}
let $g\in L^{(p^\ast)^\prime}(U,\mu)$,  
$(p^\ast)^\prime$ denotes the conjuate exponent of $p^\ast$,
and $u_j(\cdot,\om)$ be the minimizer of
$$
\min\left\{\FFepsi(u,\om)-\int_Ugu\,d\mu:\,u\in u_0+\WspoG(U,\mu)\right\},
$$
then $(u_j)$ converges weakly in $\Wsp(U,\mu)$ and $\mmu$ a.s. in $\Om$
to the minimizer of
$$
\min\left\{\FFpsi(u)-\int_Ugu\,d\mu:\,u\in u_0+\WspoG(U,\mu)\right\}.
$$
In addition, if $U=\R^N_+$ and  $g\in L^{(p^\ast)^\prime}(U,\mu)$
the minimizer $u_j(\cdot,\om)$ of
$$
\min\left\{\FFepsi(u,\om)-\int_{\R^N_+}gu\,d\mu:\,u\in\Wspo(\R^N_+,\mu)\right\},
$$
converges locally weakly in $\Wsp(\R^N_+,\mu)$ and $\mmu$ a.s. in $\Om$
to the minimizer of
$$
\min\left\{\FFpsi(u)-\int_{\R^N_+}gu\,d\mu:\,u\in\Wspo(\R^N_+,\mu)\right\}.
$$
\end{corollary}
\begin{remark}
Theorem~\ref{main} recovers the results established in \cite{Caf-Mel2}
for $p=2$. Indeed, in the statement there $N\geq 2$, $a\in(-1,1)$
and thus the compatibility condition between $a$ and $p$ is satisfied.
The results contained in \cite{Caf-Mel1} can also be inferred by the
method below (see Section~\ref{genztn}).
\end{remark}

\begin{remark}\label{uncont1}
In case $U$ has finite measure 
(Hp 5) is unnecessary since the constant function 
$\sup_{\UU}\psi$ satisfies it. In general, (Hp 5) suffices to ensure that 
$\Gamma\hbox{-}\liminf\FFjpsi$ is finite in some point, i.e. on $f$.
Actually, from Propositions~\ref{lb}, \ref{ub} below we get
$\Gamma\hbox{-}\lim_j\FFjpsi(f)=\FFpsi(f)$.
\end{remark}

\begin{remark}\label{uncont}
In \cite{Caf-Mel2} the obstacle function $\psi$ is taken to be 
defined on the whole of $U$ and to be $C^{1,1}(U)$, which clearly
implies $\sup_{\UU}\psi(\cdot,0)<+\infty$ if $U$ is bounded. 
The latter condition
is guaranteed also if $\UU=\R^{N-1}$ since the $\Gamma$-limit is finite 
in some point (see Remark~\ref{uncont1}).
Indeed, in such a case it follows that $\psi(\cdot,0)\vee 0\in L^p(\R^{N-1})$.
More generally, this holds whenever $\UU$ is not quasibounded.
\end{remark} 

\begin{remark}\label{trivial}
 The restrictions $a>-1$ and $p> 1+a$ avoid trivial results.
Indeed, if $a\leq-1$ or $p\leq 1+a$ then $\Wsp(U,\mu)\equiv\Wspo(U,\mu)$
(see \cite[Proposition 9.10]{Kuf}),
and the compatibility condition in (Hp 5) leads to $\psi\leq 0$.
Hence, no finite penalization term would appear in the homogenization limit.
\end{remark}

\section{Sobolev Spaces with $\mathcal{A}_p$ weights}\label{sobolev}

\subsection{Generalities}

We recall that a function $w\in L^1_{\mathrm{loc}}(\R^N,(0,+\infty))$
is called a \emph{Muckenhoupt $p$-weight} for $p\in(1,+\infty)$,
in short $u\in\mathcal{A}_p(\R^N)$, if
$w\in L^1_{\mathrm{loc}}(\R^N,(0,+\infty))$ and
\begin{equation}
  \label{eq:muck}
\sup_{r>0,\, z\in\R^N}\left(r^{-N}\int_{B_r(z)}w(x)\,dx\right)
\left(r^{-N}\int_{B_r(z)}w^{1/(1-p)}(x)dx\right)^{p-1}<+\infty.
\end{equation}
In the sequel we will consider only weight functions of the form
$w(x)=|\xn|^a$, with $a\in(-1,+\infty)$ and $p>(1+a)\vee 1$
(see the Appendix~\ref{Muck} and Remark~\ref{trivial}).
Then we define the Radon measure $\mu$ on $\R^N$ by $\mu:=wd\L^N$.
Take note that $\L^N\ll\mu$ and $\mu\ll\L^N$.
If $A\subseteq\R^N$ is an open set the space $H^{1,p}(A,\mu)$ defined
as the closure of $C^\infty(A)$ in $L^p(A,\mu)$ under the norm
$$
\|\varphi\|_{H^{1,p}(A,\mu)}=\|\varphi\|_{L^p(A,\mu)}+
\|\nabla\varphi\|_{(L^p(A,\mu))^N}
$$
shares several properties with the usual unweighted case.
In particular, Meyers and Serrin's $H=W$ property holds (see \cite{K2}).
We will give precise references for those properties employed in
the sequel in the respective places.
We will mainly refer to the book \cite{HKM}, and to \cite{Kuf}
when the general theory of weighted Sobolev spaces is concerned.
Hereafter we quote explicitely only those results which will be
repeatedly used in the proofs below.
\begin{lemma}\label{poincare}
  Let $A\subseteq\R^N$ be a connected bounded open set, 
then for any $E\subseteq A$ with $\L^N(E)>0$ there exists 
a constant $c=c(A,E,N,p,\mu)>0$ such that
%
\begin{equation}
  \label{eq:poincare}
  \int_{rA}|u-u_{rE}|^pd\mu\leq c r^p\int_{rA}|\nabla u|^pd\mu
\end{equation}
for any $r>0$ and $u\in \Wsp(rA,\mu)$, where $u_{rE}=\fint_{rE}u\,d\mu$.
\end{lemma}
The (scaled) Poincar\'e inequality stated above
can be inferred by the usual proof by contradiction in case $r=1$
and a simple scaling argument (see \cite[Theorem 1.31]{HKM} for
weak compactness results in weighted Sobolev spaces).
Let us then establish a weighted Sobolev-Gagliardo-Nirenberg 
inequality. 
\begin{lemma}\label{Gagliardo}
Let $a\in(-1,+\infty)$, $p\in((1+a)\vee 1,N+a)$. 
There exists a constant $c=c(N,p,\mu)>0$ such that
\begin{equation}
  \label{eq:gagliardo}
\|u\|_{L^{p^\ast}(\R^N,\mu)}\leq c \|\nabla u\|_{(L^p(\R^N,\mu))^N}
\end{equation}
for all $u\in\Kpmu$, where $p^\ast=(N+a)p/(N+a-p)$ (see \eqref{kpmu}).
\end{lemma}
\begin{proof}
 Let us first notice that the measure $\mu$ is \emph{$p$-admissible} 
according to \cite[Chapters 1,5]{HKM}. By \cite[Theorem 15.21]{HKM} 
there exist constants $\chi>1$ and $c=c(N,p,\mu)>0$ such that
\begin{equation}
  \label{eq:gagliardo0}
\left(\fint_{B_r}|u|^{\chi p}d\mu\right)^{1/(\chi p)}\leq c\, r
\left(\fint_{B_r}|\nabla u|^pd\mu\right)^{1/p}
\end{equation}
for all $r>0$ and $u\in C^\infty_0(B_r)$. 
Being the measure $\mu=|\xn|^ad\L^N$ $(N+a)$-homogeneous 
a scaling argument shows that $\chi=p^\ast/p$.
Thus, \eqref{eq:gagliardo0} rewrites as \eqref{eq:gagliardo}
for all $r>0$ and $u\in C^\infty_0(B_r)$.

The equality $\Wsp(\R^N,\mu)=\Wspo(\R^N,\mu)$ 
(see \cite[Theorem 1.27]{HKM}) and \eqref{eq:gagliardo0} 
justify \eqref{eq:gagliardo} for Sobolev maps by a 
density argument.
Eventually, given $u\in\Kpmu$ let $\varphi_n$ be a cut-off 
function between $B_{n}$ and $B_{2n}$ with 
$\|\nabla\varphi_n\|^p_{(L^\infty(\R^N))^N}\leq 2/n^p$, we claim that  
$u_n=\varphi_nu\in\Wsp(\R^N,\mu)$ converges strongly to $u$ in 
$L^{p^\ast}(\R^N,\mu)$ and $\nabla u_n$ converges strongly to 
$\nabla u$ in $(L^p(\R^N,\mu))^N$. 
Indeed, we have
$$
\int_{\R^N}|u_n-u|^{p^\ast}d\mu\leq\int_{\R^N\setminus B_n}|u|^{p^\ast}d\mu
$$
and by H\"older's inequality
\begin{eqnarray*}
\lefteqn{\int_{\R^N}|\nabla(u_n-u)|^pd\mu\leq 
2^{p-1}\int_{\R^N\setminus B_n}|\nabla u|^pd\mu
+\frac{2^p}{n^p}\int_{B_{2n}\setminus B_n}|u|^pd\mu}\\
&&\leq 2^{p-1}\int_{\R^N\setminus B_n}|\nabla u|^pd\mu
+2^p(\mu(B_2\setminus B_1))^{p/(N+a)}
\left(\int_{B_{2n}\setminus B_n}|u|^{p^\ast}d\mu\right)^{p/p^\ast}
\end{eqnarray*}
and so the conclusion follows.
\end{proof}
Finally, we recall a trace result in the weighted setting.
\begin{theorem}[Theorem 9.14 \cite{Kuf}, Sec. 10.1 \cite{Nik}]\label{thmtracce}
Let $A\subseteq\R^N$ be a Lipschitz bounded open set, 
if $a\in(-1,+\infty)$ and $p>1+a$ there exists a compact operator
$\mathrm{Tr}:\Wsp(A,\mu)\to L^p(\partial_NA)$ such that $\mathrm{Tr}(u)=u$
for every $u\in C^\infty(\overline{A})$.
\end{theorem}
In the rest of the paper 
to denote the trace of a function $u\in \Wsp(A,\mu)$ on 
$\partial_NA$ we use the more appealing notation $u(\cdot,0)$.

\subsection{Variational $(p,\mu)$-capacities}
We recall the notion of variational $(p,\mu)$-capacity
(see \cite[Chapter 2]{HKM}): Given any open set $A\subseteq\R^N$
and any set $E\subseteq\R^N$ define
$$
\capmu(E,A):=\inf_{\{A^\prime \text{ open }:\,A^\prime \supseteq E\}}
\inf\left\{\int_{A}|\nabla u|^pd\mu:\,
u\in\Wspo(A,\mu),\,u\geq 1\,\L^N \text{ a.e. on } A^\prime\right\},
$$
with the usual convention $\inf\emptyset=+\infty$.
In case $A=\R^N$, $N\geq 2$, we drop the dependence
on $A$ and write only $\capmu(E)$.

Recall that a property holds $\capmu$ q.e. if it holds up to
a set of $\capmu$ zero. In particular, any function $u$ in
$\Wsp(A,\mu)$ has a \emph{precise representative} $\tilde{u}$
defined $\capmu$ q.e. (see \cite[Chapter 4]{HKM} and \cite{K1}).
By means of this result the following formula holds
(see \cite[Corollary 4.13]{HKM} and the subsequent comments)
\begin{equation}
  \label{eq:capaltern}
\capmu(E,A)=\inf\left\{\int_{A}|\nabla u|^pd\mu:\,
u\in\Wspo(A,\mu),\,\tilde{u}\geq 1\,\capmu \text{ q.e. on } E\right\}.
\end{equation}
Thanks to \eqref{eq:capaltern} it is easy to show that   
if $A$ is bounded the minimum problem for the capacity 
has a unique minimizer $u^{E,A}$, called the 
{\sl $(p,\mu)$-capacitary potential} of $E$ in $A$. 
Instead, in case $A=\R^N$ the minimizer might not exist. 
The minimum problem has to be relaxed, so that it has 
a (unique) solution, denoted by $u^E$, in the space 
$\Kpmu$ by Lemma~\ref{Gagliardo}.

Simple truncation arguments imply that
$0\leq u^E\leq 1$ $\L^N$ a.e. on $\R^N$, and
for every $\lambda>0$ we get by scaling
\begin{equation}
  \label{eq:cap-scaling}
\capmu(\lambda E,\lambda A)=\lambda^{N-p+a}\capmu(E,A).
\end{equation}
For this reason we will restrict ourselves to the range $p<N+a$ to be 
sure that points have zero capacity (see for instance
\cite[Theorem 2.19]{HKM}).

If $A$ and $E$ are simmetric with respect to the hyperplane
$\xn=0$ then the $(p,\mu)$-capacitary potential of $E$ 
in $A$, $u^{E,A}$, enjoys the same simmetry and in 
addition it satisfies
\begin{equation}
  \label{eq:cap-sim}
  \int_{A\cap\R^N_+}|\nabla u^{E,A}|^pd\mu=\frac 12\capmu(E,A).
\end{equation}
Moreover, $\capmu(E+z,A+z)=\capmu(E,A)$ if $z\in\R^{N-1}\times\{0\}$, 
being $\mu$ unaffected by horizontal translations.

Some further properties are needed. The results below are elementary,
but since we have found no explicit reference in literature we
prefer to give full proofs.

First we show that set inclusion induces a partial ordering
among capacitary potentials.
\begin{proposition}\label{cappotential}
 Assume $E\subseteq F$, then $u^E\leq u^F$ $\L^N$ a.e. in $\R^N$.
\end{proposition}
\begin{proof}
Assume by contradiction that $\L^N(\{u^F<u^E\})>0$, then the test-function
$\varphi=(u^E-u^F)\vee 0\in \Wspo(\R^N,\mu)$ is not identically $0$.
Notice that
\begin{eqnarray}\label{eq:grad}
\nabla\varphi=
  \begin{cases}
  \nabla(u^E-u^F) & \L^N \text{ a.e. in } \{u^F<u^E\} \\
0 & \L^N \text{ a.e. in } \{u^E\leq u^F\} 
  \end{cases}
\end{eqnarray}
(see \cite[Theorem 1.20]{HKM}).
By exploiting the strict minimality of $u^F$ for the capacitary
problem related to $F$, and by comparing its energy with that of
$u^F+\varphi$, \eqref{eq:grad} entails
\begin{equation}
  \label{eq:cap}
\int_{\{u^F<u^E\}}|\nabla u^F|^pd\mu<\int_{\{u^F<u^E\}}|\nabla u^E|^pd\mu.
\end{equation}
Let us now define $w=u^E\wedge u^F$, then $w$ is admissible for the
capacitary problem related to $E$, and by computing its energy
we infer from \eqref{eq:cap}
\begin{eqnarray*}
\int_{\R^N}|\nabla w|^pd\mu=
\int_{\{u^E\leq u^F\}}|\nabla u^E|^pd\mu+
\int_{\{u^F<u^E\}}|\nabla u^F|^pd\mu
<\int_{\R^N}|\nabla u^E|^pd\mu,
\end{eqnarray*}
which is clearly a contradiction.
\end{proof}
In turn, Proposition~\ref{cappotential} yields uniform convergence
of the relative capacities to the global one for sets contained in a
bounded open given one.
In doing that we exploit De Giorgi's slicing-averaging method to refine
the cut-off argument contained in Lemma~\ref{Gagliardo}.
\begin{proposition}\label{cap-prop}
For any bounded set $E\subset\R^N$ we have
\begin{equation}
  \label{eq:rel-cap-conv}
  \lim_n\capmu(E,B_n)=\inf_n\capmu(E,B_n)=\capmu(E).
\end{equation}
Furthermore, given a bounded open set $A\subseteq\R^N$, then
  \begin{equation}
    \label{eq:rel-cap-unif}
  \lim_n\sup_{\{E:\,E\subseteq A\}}|\capmu(E,B_n)-\capmu(E)|=0.
  \end{equation}
\end{proposition}
\begin{proof}
 Assume $E\subset\subset B_m$, and let $\varphi_n^k$ be a 
cut-off function between $B_{nk}$ and $B_{n(k+1)}$, $n,r\in\N$ 
with $n\geq r\geq m$ and $k\in\{1,\ldots,r-1\}$, 
such that $\|\nabla\varphi_n^k\|^p_{(L^\infty(\R^N))^N}\leq 2/n^p$.
Thus, it follows for every such $k$
\begin{eqnarray*}
\lefteqn{ \capmu(E,B_{n(k+1)})\leq\int_{\R^N}|\nabla(\varphi_n^k u^E)|^pd\mu}
\\&&\leq\int_{B_{nk}}|\nabla u^E|^pd\mu+
2^{p-1}\int_{B_{n(k+1)}\setminus B_{nk}}|\nabla u^E|^pd\mu+
\frac{2^p}{n^p}\int_{B_{n(k+1)}\setminus B_{nk}}|u^E|^pd\mu.
\end{eqnarray*}
Hence, by taking into account that $(\capmu(E,B_i))_{i\in\N}$ is a
decreasing sequence and by summing-up on $k\in\{1,\ldots,r-1\}$ and 
averaging we infer
$$
\capmu(E,B_{nr})\leq\left(1+\frac{2^p}r\right)\capmu(E)
+\frac{2^p}{rn^p}\int_{B_{nr}\setminus B_{n}}|u^E|^pd\mu.
$$
Since $u^E\in L^{p^\ast}(\R^N,\mu)$ by Lemma~\ref{Gagliardo},  
H\"older's inequality and the $(N+a)$-homogeneity of $\mu$ yield
\begin{equation}\label{eq:cap-rel-stima}
\capmu(E,B_{nr})\leq\left(1+\frac{2^p}r\right)\capmu(E) 
+2^p\mu(B_1)\,r^{p-1}
\left(\int_{B_{nr}\setminus B_{n}}|u^E|^{p^\ast}d\mu\right)^{p/p^\ast}.
\end{equation}
In turn, by passing to the limit first as $n\to+\infty$ and 
then as $r\to+\infty$ the latter estimate implies \eqref{eq:rel-cap-conv} 
being $(\capmu(E,B_i))_{i\in\N}$ decreasing and bounded from below
by $\capmu(E)$.

Eventually, to get \eqref{eq:rel-cap-unif} notice that
with fixed a bounded open set $A$, $A\subset\subset B_m$,
for every $E\subseteq A$ we have $\capmu(E)\leq\capmu(A)$ and
$0\leq u^E\leq u^A$ by Proposition~\ref{cappotential}, then
\eqref{eq:cap-rel-stima} yields
$$
0\leq \capmu(E,B_{nr})-\capmu(E)\leq\frac{2^p}r\capmu(A)
+2^p\mu(B_1)\,r^{p-1}
\left(\int_{B_{nr}\setminus B_{n}}|u^A|^{p^\ast}d\mu\right)^{p/p^\ast}.
$$
By taking into account that $(\capmu(E,B_i))_{i\in\N}$ is decreasing 
the uniform convergence is established .
\end{proof}

\section{A weighted Ergodic Theorem}\label{weighted-ergodic}

In this section we prove a weighted version of the ergodic theorem relevant
in our analysis. We adopt the notation of (Hp 2) and introduce
some new. First, take note that $\gamma(\ii,\cdot)\in  L^\infty(\Om,\mmu)$
for every $\ii\in\Z^{N-1}$, and that the stationarity assumption
\eqref{eq:ergogamma} on the $\tau_{\ii}$'s yields
$\mathbb{E}[\gamma(\ii,\cdot)]=\mathbb{E}[\gamma(\kk,\cdot)]$
for every  $\ii,\kk\in\Z^{N-1}$, where
$$
\mathbb{E}[\gamma(\ii,\cdot)]:=\int_\Om\gamma(\ii,\om)\,d\mmu(\om).
$$
The common value is denoted simply by $\mathbb{E}[\gamma]$.
For every $\ii\in\Z^{N-1}$ the operator
$T_{\ii}:L^\infty(\Om,\mmu)\to L^\infty(\Om,\mmu)$ is defined by
$T_{\ii}(f)=f\circ\tau_{\ii}$.
By the stationarity assumption \eqref{eq:ergogamma} it is then
easy to check that $\mathscr{S}=\left\{T_{\ii}\right\}_{\ii\in\Z^{N-1}}$
is a multiparameter semigroup generated by the commuting isometries
$T_{{\tt e}_r}$ for $r\in\{1,\ldots,N-1\}$, being
$\{{\tt e}_1,\ldots,{\tt e}_{N-1}\}$ the canonical basis of $\R^{N-1}$.
\begin{theorem}\label{ergo-w}
Let $\gamma$ be a process satisfying (Hp 2), then $\mmu$ a.s. in $\Om$
\begin{equation}
  \label{eq:Birkhoff0}
\lim_j\frac{1}{\#\Ieps(V)}\sum_{\ii\in\Ieps(V)}\gammaio=\mathbb{E}[\gamma],
\end{equation}
and
\begin{equation}
  \label{eq:Birkhoff1}
\Psi_j(x,\om):=\sum_{\ii\in\Ieps(V)}\gammaio\chi_{\Qijo}(x)\to\mathbb{E}[\gamma]
\quad  \mathrm{weak}^\ast\, L^\infty(V).
\end{equation}
for every bounded open set $V\subset\R^{N-1}$ with $\LN(\partial V)=0$.
\end{theorem}
\begin{proof}
With fixed a set $V$ as in the statement above define
$$
A_{j}(f)=\sum_{\ii\in\Z^{N-1}}\alpha_{j,\ii}T_\ii(f),
$$
where for every $j\in\N$ and $\ii\in\Z^{N-1}$ we set
$\alpha_{j,\ii}=(\#\Ieps(V))^{-1}\chi_{\Ieps(V)}(\ii)$.
We claim that $(A_j)_{j\in\N}$ is an ergodic $\mathscr{S}$-net
according to \cite[p.75]{Kr}, i.e.
\begin{itemize}
\item[{(E1)}] each $A_j$ is a linear operator on $L^\infty(\Om,\mmu)$,

\item[{(E2)}] $A_j(f)\in\overline{\mathrm{co}}\,\mathscr{S}(f)$
for each $f\in L^\infty(\Om,\mmu)$ and all $j\in\N$,

\item[{(E3)}] the $A_j$'s are equi-continuous, and

\item[{(E4)}]  for each $f\in L^\infty(\Om,\mmu)$ and $\ii\in\Z^{N-1}$
$$
\lim_j(A_j(T_{\ii}(f))-A_j(f))=\lim_j(T_{\ii}(A_j(f))-A_j(f))=0.
$$
\end{itemize}
Clearly, (E1) is satisfied.
For what (E2) and (E3) are concerned it is enough to notice that
$\sum_{\ii\in\Z^{N-1}}\alpha_{j,\ii}=1$ for every $j\in\N$.
Moreover, for $j$ sufficiently big it holds
\begin{eqnarray*}
\lefteqn{\sum_{\ii\in\Z^{N-1}}\sum_{\{\kk\in\Z^{N-1}:\,|\kk|=1\}}
|\alpha_{j,\ii}-\alpha_{j,\ii+\kk}|}\\&&\leq N
\frac{\#\{\ii\in\Z^{N-1}:\,Q^{\ii}_j\cap\partial V\neq\emptyset\}}{\#\Ieps(V)}
\leq N\frac{\LN((\partial V)_{\sqrt{N}\eps})}
{\LN(V\setminus(\partial V)_{2\sqrt{N}\eps})}.
\end{eqnarray*}
In turn the latter estimate implies
$$
\limsup_j\left(|\alpha_{j,0}|+\sum_{\ii\in\Z^{N-1}}\sum_{\{\kk\in\Z^{N-1}:\,|\kk|=1\}}
|\alpha_{j,\ii}-\alpha_{j,\ii+\kk}|
\right)=0
$$
and thus (E4) is satisfied, too.
By Eberlein's Theorem (see \cite[Theorem 1.5 p.76]{Kr}) we have that
$A_j(f)\to \bar{f}$ in $L^\infty(\Om,\mmu)$ for all $f\in L^\infty(\Om,\mmu)$,
where $\bar{f}\in\{g\in\overline{\mathrm{co}}\,\mathscr{S}(f):\,
T_{\ii}(g)=g\,\, \forall \ii\in\Z^{N-1}\}$.
The ergodicity assumption on the $\tau_{\ii}$'s implies that
$\bar{f}$ is constant  $\mmu$ a.s. in $\Om$, and since
$\sum_{\ii\in\Z^{N-1}}\alpha_{j,\ii}=1$ for every $j\in\N$,
the convergence $A_j(f)\to\bar{f}$ in $L^\infty(\Om,\mmu)$
implies $\bar{f}=\mathbb E[f]$.

To deduce \eqref{eq:Birkhoff0} apply the result above to
$\gamma(\underline{0},\cdot)$ and notice that
$$
A_j(\gamma(\underline{0},\om))=
\frac{1}{\#\Ieps(V)}\sum_{\ii\in\Ieps(V)}\gamma(\ii,\om)
$$
since $\gammaio=\gamma(\underline{0},\tau_{\ii}(\om))$
for all $\ii\in\Z^{N-1}$ and $\om\in\Om$ by \eqref{eq:ergogamma}.

Eventually, in order to prove \eqref{eq:Birkhoff1}
consider the family $\mathscr{Q}$ of all open cubes in $\R^{N-1}$
with sides parallel to the coordinate axes, and with center and vertices
having rational coordinates.
To show the claimed weak$^\ast$ convergence it suffices to check that
$\lim_j\int_\Om\Psi_j(x,\om)\chi_Q(x)\,d\LN=\LN(Q)\mathbb{E}[\gamma]$
for any $Q\in\mathscr{Q}$ with $Q\subseteq V$.
We have
\begin{eqnarray*}
\lefteqn{\left|\int_{Q}(\Psi_j(x,\om)-
\mathbb{E}[\gamma])d\LN\right|}\\&&
\leq\left|\epsN\sum_{\ii\in\Ieps(Q)}\gammaio
-\LN(Q)\mathbb{E}[\gamma]\right|+2\gammab
\LN\left(Q\setminus\cup_{\ii\in\Ieps(Q)}\Qijo\right),
\end{eqnarray*}
and thus \eqref{eq:Birkhoff0} and the denumerability of $\mathscr{Q}$
yield that the rhs above is infinitesimal $\mmu$ a.s. in $\Om$.
\end{proof}

\begin{remark}\label{only-stationary}
Even dropping the ergodicity assumption, 
conclusions similar to 
those in Theorem~\ref{ergo-w} still hold true. 
Indeed, by arguing as in the proof above integrating and  
exploiting the stationarity of $\gamma$, the limit 
$\bar{\gamma}(\underline{0},\cdot)$ of the sequence 
$(A_j(\gamma(\underline{0},\cdot)))$ turns out 
to be characterized as the unique 
function in $L^\infty(\Om,\mmu)$ satisfying 
$$
\int_I\gamma(\underline{0},\om)d\mmu=
\int_I\bar{\gamma}(\underline{0},\om)d\mmu
$$
for every set $I\in\mathscr{P}$ invariant w.r.to the
$\tau_\ii$'s. Thus, if 
$\mathscr{I}$ denotes the $\sigma$-subalgebra  of $\mathscr{P}$ 
of the invariant sets of the $\tau_\ii$'s, 
$\bar{\gamma}(\underline{0},\cdot)$ is 
the conditional expectation of $\gamma(\underline{0},\cdot)$
relative to $\mathscr{I}$, denoted by $\mathbb{E}[\gamma,\mathscr{I}]$.
Statement \eqref{eq:Birkhoff1} then follows analogously. 

\end{remark}

\section{Proof of the main result}\label{mainresult}

Throughout the section the open set $U\subseteq\R^N_+$ will be
fixed. Thus, for the sake of simplicity we denote $\Ieps:=\Ieps(\UU)$.
Furthermore, $(V_n)_{n\in\N}$ will always denote a sequence of bounded
open subsets of $\UU$ with Lipschitz boundary such that $\UU=\cup_nV_n$ 
and $V_n\subset\subset V_{n+1}$.

The set $\Om^\prime$ mentioned in Theorem~\ref{main} is defined
as any subset of $\Om$ of full probability for which \eqref{eq:Birkhoff0}
and \eqref{eq:Birkhoff1} hold true for $V_n$ for every $n\in\N$.

In some computations we find inequalities involving constants
depending on $U$, $N$, $p$, $\mu$ etc... but are always independent
from the indexing parameter $j$. Since it is not essential to
distinguish from one
specific constant to another, we indicate all of them by the same letter $c$,
leaving understood that $c$ may change from one inequality to another.

Below we prove a joining lemma on varying boundary domains for weighted
Sobolev type energies. The argument follows closely that by \cite{ANB}
in the unweighted case for the periodic homogenization on perforated
open sets.
\begin{lemma}\label{joining}
Let $(u_j)$ be converging to $u$ in $L^p(U,\mu)$ for which
$\sup_j\|u_j\|_{\Wsp(U,\mu)}<+\infty$.

Let $k\in\N$ and $\om\in\Om$ be fixed, then for all $\ii\in \Ieps$
there exists $h_\ii\in\{1,\ldots,k\}$ such that,
having set
\begin{eqnarray*}
\Bijh:=\{x\in U:\, |x-\xije|<2^{-h}\eps\},\quad
\Cijh:=\Bijh\setminus B_j^{i,h+1},
\end{eqnarray*}
there exists a sequence $(\vj)$ converging to $u$ in
$L^p(U,\mu)$ and such that for every $j\in\N$
\begin{equation}
  \label{eq:equal1}
  \vj\equiv u_j\, \text{ on }\, U\setminus\cup_{\ii\in \Ieps}\Cij,
\end{equation}
\begin{equation}
  \label{eq:equal2}
  \vj(x)\equiv (u_j)_{\Cij}\, \text{ if }\, |x-\xije|=\frac 342^{-h_\ii}\eps,\,
x\in U,
\end{equation}
\begin{equation}
  \label{eq:errorest}
\left|\int_{A}|\nabla\vj|^pd\mu-
\int_{A}|\nabla u_j|^pd\mu\right|\leq
\frac ck\int_{\cup_{\ii\in\IIeps(A)}\Qij\times(0,\eps)}|\nabla u_j|^pd\mu
\end{equation}
for some positive constant $c$ independent from $j$ and $k$, 
and for all open sets $A\subseteq U$ where
$$
\IIeps(A):=\{\ii\in\Ieps:\,\Qijo\cap A\neq\emptyset\}.
$$
Furthermore, the functions $\zeta_j:=\sum_{\ii\in \Ieps}(u_j)_{\Cij}\chi_{\Qijo}$
converge to $u$ in $L^p_{\mathrm{loc}}(\UU)$.
\end{lemma}
\begin{proof} For all $j\in\N$, $\ii\in\Ieps$ and $1\leq h\leq k$
denote by $\varphi^{\ii,h}_j$ a cut-off function between
$\Sijh:=\{x\in U:\, |x-\xije|=\frac 342^{-h}\eps\}$ and
$U\setminus\Cijh$, with $\|\nabla\varphi^{\ii,h}_j\|_{(L^\infty(\R^N))^N}
\leq 2^{h+2}\eps^{-1}$. Then define
\begin{eqnarray*}
v^{\ii,h}_j:=
\begin{cases}
\varphi^{\ii,h}_j(u_j)_{C^{\ii,h}_j}+(1-\varphi^{\ii,h}_j)u_j &
\text{ on } \Cijh,\,\ii\in\Ieps\\
u_j & \text{ otherwise on } U.
\end{cases}
\end{eqnarray*}
Being $\mathrm{Lip}(\varphi^{\ii,h}_j)2^{-h-2}\eps\leq 1$ we infer
\begin{eqnarray*}
\|\nabla v^{\ii,h}_j\|_{(L^p(\Cijh,\mu))^N}^p\leq
\|\nabla u_j\|_{(L^p(\Cijh,\mu))^N}^p+
\left(\frac{2^{h+2}}{\eps}\right)^p
\int_{\Cijh}|u_j-(u_j)_{\Cijh}|^pd\mu,
\end{eqnarray*}
and thus by taking into account Lemma~\ref{poincare} we get
\begin{equation}\label{stima001}
\|\nabla v^{\ii,h}_j\|_{(L^p(\Cijh,\mu))^N}^p\leq
c\|\nabla u_j\|_{(L^p(\Cijh,\mu))^N}^p,
\end{equation}
for some positive constant $c$ depending only on $N$, $p$ and $\mu$.
Indeed, the ratio between the outer and inner radii of $\Cijh$ is
equals $2$ for every $\ii,j,h$.

By summing up and averaging in $h$, being the $\Cijh$ disjoint,
we find $h_\ii\in\{1,\ldots,k\}$ such that
\begin{equation}\label{stima002}
\|\nabla u_j\|_{(L^p(C_j^{\ii,h_\ii},\mu))^N}^p\leq
\frac ck\|\nabla u_j\|_{(L^p(\Qijo\times(0,\eps),\mu))^N}^p.
\end{equation}
Define $\vj=v^{\ii,h_\ii}_j$ on $\cup_{\Ieps}C_j^{\ii,h_\ii}$ and $\vj=u_j$
otherwise, then \eqref{eq:equal1}, \eqref{eq:equal2} are satisfied by
construction, and \eqref{eq:errorest} follows easily from 
\eqref{stima001} and \eqref{stima002}.

To prove that $(\vj)$ converges to $u$ in $L^p(U,\mu)$ we use again
Lemma~\ref{poincare}. Indeed, by the very definition of $v_j$ we have
\begin{eqnarray*}
\lefteqn{\|u_j-\vj\|_{L^p(U,\mu)}=
\sum_{\ii\in \Ieps}\|u_j-\vj\|_{L^p(C_j^{\ii,h_\ii},\mu)}}\\
&&\leq\sum_{\ii\in \Ieps}\|u_j-(u_j)_{C_j^{\ii,h_\ii}}\|_{L^p(C_j^{\ii,h_\ii},\mu)}\leq c
\sum_{\ii\in \Ieps}\frac{\eps}{2^{h_i}}\|\nabla u_j\|_{(L^p(C_j^{\ii,h_\ii},\mu))^N}
\leq c\eps\|\nabla u_j\|_{(L^p(U,\mu))^N}.
\end{eqnarray*}

Eventually, let us show the convergence of $(\zeta_j)$ to $u$ in 
$L^p_{\mathrm{loc}}(\UU)$.
The (local) compactness of the trace operator (see Theorem~\ref{thmtracce})
and the very definition of $v_j$ entail for any open set 
$V\subset\subset\UU$ 
\begin{equation}
  \label{eq:uno}
\limsup_j\|\zeta_j-u\|_{L^p(V)}^p=\limsup_j\|\zeta_j-u_j\|_{L^p(V)}^p\leq
\limsup_j\sum_{\ii\in \Ieps}\|u_j-(u_j)_{\Cij}\|_{L^p(\Qijo)}^p.
\end{equation}
An elementary scaling argument and the Trace theorem~\ref{thmtracce} yield
\begin{eqnarray}\label{eq:tre}
\lefteqn{\|u_j-(u_j)_{\Cij}\|_{L^p(\Qijo)}^p}\notag\\&&
\leq c\eps^{-1-a}\left(
\|u_j-(u_j)_{\Cij}\|_{L^p(\Qijo\times(0,\eps),\mu)}^p+
\eps^p\|\nabla u_j\|_{(L^p(\Qijo\times(0,\eps),\mu))^N}^p
\right)
\end{eqnarray}
for some positive constant $c$ depending only on $N$, $p$ and $\mu$.
Since the scaled Poincar\'e inequality \eqref{eq:poincare} entails
\begin{equation}
  \label{eq:quattro}
\|u_j-(u_j)_{\Cij}\|_{L^p(\Qijo\times(0,\eps),\mu)}^p\leq c
\eps^p\|\nabla u_j\|_{(L^p(\Qijo\times(0,\eps),\mu))^N}^p,
\end{equation}
the thesis then follows by collecting \eqref{eq:uno}-\eqref{eq:quattro}
being $p>1+a$.
\end{proof}
We are now ready to prove the lower bound inequality.

\begin{proposition}\label{lb}
For all $\om\in\Om^\prime$ and $u\in L^p(U,\mu)$
\begin{equation}
  \label{eq:lb}
\FFpsi(u)\leq\Gamma\hbox{-}\liminf_j\FFjpsi(u,\om).
\end{equation}
\end{proposition}
\begin{proof}
We may assume $\Lambda\in(0,+\infty)$, the estimate being trivial
if $\Lambda$ equals $0$, while if $\Lambda=+\infty$ it
can be inferred by a simple comparison argument with the case
$\Lambda$ finite.

We use the notation introduced in Lemma~\ref{joining}, and further set
\begin{equation}
  \label{eq:bji}
  \Bji:=\left\{x\in\R^N:\, |x-\xije|<\frac 342^{-h_\ii}\eps\right\},
\end{equation}
for all $\ii\in\Ieps$ (recall that $\Ieps=\Ieps(\UU)$).

With fixed $u_j\to u$ in $L^p(U,\mu)$ with $\sup_j\FFjpsi(u_j,\om)<+\infty$
define the function
\begin{eqnarray*}
\xi_j(x):=\begin{cases}
(u_j)_{C^{\ii,h_\ii}_j} & \text{ on } \Bji\cap U,\,\ii\in\Ieps\\
\vj(x) & \text{ otherwise on } U,
\end{cases}
\end{eqnarray*}
where $(\vj)$ is the sequence provided by Lemma~\ref{joining}.
It is easy to check that $\xi_j\to u$ in $L^p(U,\mu)$ and
$\sup_j\|\xi_j\|_{\Wsp(U,\mu)}<+\infty$.
By taking into account \eqref{eq:errorest} and by splitting the energy
contribution of $\vj$ far from and close to the obstacles yields
\begin{eqnarray}\label{parziale}
\lefteqn{\left(1+\frac ck\right)\liminf_j\FFjpsi(u_j,\om)\geq
\liminf_j\FFjpsi(\vj,\om)\notag}\\&&\geq
\liminf_j\|\nabla \vj\|^p_{(L^p(U\setminus\cup_{\Ieps}\Bji,\mu))^N}+
\liminf_j\sum_{\ii\in \Ieps}\|\nabla \vj\|^p_{(L^p(\Bji\cap U,\mu))^N}\notag\\&&
=\liminf_j\|\nabla \xi_j\|^p_{(L^p(U,\mu))^N}+
\liminf_j\sum_{\ii\in \Ieps}\|\nabla \vj\|^p_{(L^p(\Bji\cap U,\mu))^N} \notag\\&&
\geq\|\nabla u\|^p_{(L^p(U,\mu))^N}+
\liminf_j\sum_{\ii\in \Ieps}\|\nabla \vj\|^p_{(L^p(\Bji\cap U,\mu))^N}.
\end{eqnarray}
We claim that for all $\om\in\Om^\prime$ 
\begin{equation}
  \label{eq:capacity}
\liminf_j\sum_{\ii\in \Ieps}\|\nabla \vj\|^p_{(L^p(\Bji\cap U,\mu))^N}\geq
\frac 12\Lambda\mathbb{E}[\gamma]\int_{\UU}\Phi(\psi(\xnc)-u(\xnc,0))\,d\xnc,
\end{equation}
where $\Phi(t):=(t\vee 0)^p$. 
Given this for granted, we infer \eqref{eq:lb} 
from \eqref{eq:capacity} and by letting $k\to+\infty$ in
\eqref{parziale}.

To conclude we are left with proving \eqref{eq:capacity}. Denote by 
$\hat{U}=\mathrm{int}
\{(\xnc,x_N)\in\R^N:\,(\xnc,|x_N|)\in\overline{U}\}$, 
and extend $v_j$ to $\hat{U}$ by simmetry with respect to 
the plane $x_N=0$, i.e. 
$\hat{v}_j(\xnc,x_N):=v_j(\xnc,|x_N|)$ for $x\in\hat{U}$.
Notice that $\hat{v}_j\in\Wsp(\hat{U},\mu)$
and $\|\nabla\hat{v}_j\|_{(L^p(\Bji,\mu))^N}^p=
2\|\nabla v_j\|_{(L^p(\Bji\cap U,\mu))^N}^p$.
Thus, for every $\ii\in\Ieps$ we infer by property 
\eqref{eq:equal2} in Lemma~\ref{joining} 
\begin{eqnarray*}  
\lefteqn{\|\nabla\vj\|_{(L^p(\Bji\cap U,\mu))^N}^p=
\frac 12\|\nabla\hat{v}_j\|_{(L^p(\Bji,\mu))^N}^p}\\&&\geq\frac 12
\inf\left\{\|\nabla v\|_{(L^p(\Bji,\mu))^N}^p:\,v-(u_j)_{\Cij}\in\Wspo(\Bji,\mu),\,
\tilde{v}\geq\psi\,\,\capmu\text{ q.e. on } \Teio\right\}\notag\\&&
\geq\frac 12
\inf\left\{\|\nabla v\|_{(L^p(\R^N,\mu))^N}^p:\,v\in\Wspo(\R^N,\mu),\,
\tilde{v}\geq\psi-(u_j)_{\Cij}\,\,\capmu\text{ q.e. on } \Teio\right\}.\notag
\end{eqnarray*}
With fixed $\eta>0$ the uniform continuity of $\psi$ on the open set
$V_{n+1}\subset\subset \UU$ implies that $\psi(y)\geq\psi(\xijo)-\eta$
for every $y\in\cup_{\Ieps(V_{n+1})}\Teio$ for $j$ sufficiently big.
Thus we deduce
\begin{equation}\label{eq:finale}
\sum_{\ii\in\Ieps}\|\nabla \vj\|^p_{(L^p(\Bji,\mu))^N}
\geq\frac 12\deltaj\sum_{\ii\in \Ieps(V_{n+1})}
\gammaio\Phi(\psi(\xijo)-(u_j)_{\Cij}-\eta).
\end{equation}
In deriving the last inequality we have exploited the $p$-homogeneity
of the weighted norm, formula 
\eqref{eq:capaltern}, and the capacitary scaling assumption in (Hp 1).

To estimate the last term above define
$\psi_j:=\sum_{\ii\in \Ieps}(\psi(\xijo)-(u_j)_{\Cij})\chi_{\Qijo}$
and consider the functions $\Psi_j$ introduced in
Theorem~\ref{ergo-w} for $V=V_{n+1}$, i.e.
$$
\Psi_j(\xnc,\om)=\sum_{\ii\in\Ieps(V_{n+1})}\gammaio\chi_{\Qijo}(\xnc).
$$
Recall that by the very definition of
$\Om^\prime$ we have $\Psi_j(\cdot,\om)\to\mathbb{E}[\gamma]$
weak$^\ast$ $L^\infty(V_{n+1})$ for all $\om\in\Om^\prime$.
Being $V_n\subset\subset V_{n+1}$, \eqref{eq:finale} rewrites
for $j$ sufficiently big as
\begin{eqnarray}\label{duality}
\sum_{\ii\in\Ieps}\|\nabla \vj\|^p_{(L^p(\Bji,\mu))^N}\geq
\frac 12
\deltaj\eps^{-N+1}\int_{V_n}\Phi(\psi_j(\xnc)-\eta)\Psi_j(\xnc,\om)\,d\xnc.
\end{eqnarray}
Notice that $\psi_j\to(\psi-u)$ in $L^p(V_{n+1})$
by the continuity of $\psi$ and by Lemma~\ref{joining}.
In turn this implies $\Phi(\psi_j-\eta)\to \Phi(\psi-u-\eta)$
in $L^1(V_{n+1})$ for every $\eta>0$.
Hence, for every $k\in\N$ and $\eta>0$ we get
\begin{eqnarray*}
\liminf_j\sum_{\ii\in \Ieps}\|\nabla\vj\|^p_{(L^p(\Bji,\mu))^N}\geq\frac 12
\Lambda\mathbb{E}[\gamma]\int_{V_n}\Phi(\psi(\xnc)-u(\xnc,0)-\eta)\,d\xnc.
\end{eqnarray*}
To recover \eqref{eq:capacity} let $\eta\to 0^+$, 
and then increase $V_n$ to $\UU$. 
\end{proof}

In the next proposition we show that the lower bound derived
in Proposition~\ref{lb} is sharp.
\begin{proposition}\label{ub}
For all $\om\in\Om^\prime$ and $u\in L^p(U,\mu)$
\begin{equation}
  \label{eq:ub}
\Gamma\hbox{-}\limsup_j\FFjpsi(u,\om)\leq\FFpsi(u).
\end{equation}
\end{proposition}
\begin{proof}
  Let us show that for every $u\in L^p(\Om,\mu)$ such that
$\FFpsi(u)<+\infty$ and for every event $\om\in\Om^\prime$ we may
construct $u_j\in\Wsp(U,\mu)$ such that $u_j\to u$ in $L^p(U,\mu)$ and
\begin{equation}
  \label{eq:suppsi}
  \limsup_j\FFjpsi(u_j,\om)\leq\FFpsi(u).
\end{equation}
Take note that we may assume $\Lambda\in(0,+\infty)$. 
Indeed, if $\Lambda=0$
we may use a comparison argument with the former case to conclude.
Instead, if $\Lambda=+\infty$ by Proposition~\ref{lb} we have 
$\tilde{u}\geq\psi$ $\capmu$ q.e on $\UU$, and then we may take 
$u_j\equiv u$.

Furthermore, we may reduce to $u\in C^{0,1}\cap L^\infty(U)$, 
and $\psi\in L^\infty(\UU)$ and continuous in the relative 
interior of $\UU$ w.r.to the relative topology of $\{\xn=0\}$.

Indeed, suppose \eqref{eq:suppsi} proven under those assumptions.
The functions $\psi_k:=\psi\vee(-k)$, 
$k\in\N$, are bounded and continuous on the relative 
interior of $\UU$, $\psi_k\geq\psi_{k+1}\geq\psi$ 
and $(\psi_k)$ converges to $\psi$ pointwise. 
Denote by $\FF_j^{\psi_k}$, $\FF^{\psi_k}$ the functionals defined as $\FF_j$, 
$\FF$ in \eqref{eq:fapprox} and \eqref{eq:Glimit}, respectively, 
with $\psi$ substituted by $\psi_k$. 
Clearly, we have $\FF_j\leq\FF_j^{\psi_k}$, so that
$\Gamma\hbox{-}\limsup_j\FFjpsi(u,\om)\leq
\Gamma\hbox{-}\limsup_j\FF_j^{\psi_k}(u,\om)=\FF^{\psi_k}(u)$. Moreover,
notice that $\FF^{\psi_k}(u)\to\FFpsi(u)$ as $k\to+\infty$ being 
$u\in L^\infty(U)$.

It is easy to check that if $\FFpsi(u)<+\infty$ then
$\FFpsi(u\vee(-k)\wedge k)<+\infty$ for any $k\in\N$ with 
$k\geq\|\psi\|_{L^\infty(\UU)}$ 
(see \cite[Lemma 1.19]{HKM} for the fact that truncations 
preserve $\Wsp(U,\mu)$ regularity). 
The density of $C^{0,1}\cap L^\infty(U)$ in $\Wsp(U,\mu)$
and the lower semicontinuity of $\Gamma\hbox{-}\limsup_j\FFjpsi$
then establish \eqref{eq:suppsi} for functions in $\Wsp(U,\mu)$ once 
it has been proven for their truncations (see \cite[Theorem 1.1]{Chua} and 
\cite[Theorem 4]{K2} for extension and density results in $\R^N$, 
respectively).

Clearly, if $\psi$ is bounded we may also take the function $f$ in (Hp 5)
to be in $L^\infty(U)$ upon substituting it with 
its truncation at the levels $\pm\|\psi\|_{L^\infty(\UU)}$.

To conclude the proof we distinguish two cases according 
to whether $U$ is bounded or not.
\vskip0.25cm

\emph{Step 1: $U$ is bounded.} 
With fixed $\eta>0$ such that 
\begin{equation}
  \label{eq:regularity}
\HH\left(\{y\in\UU:\,\psi(y)-u(y,0)=\eta\}\right)=0,  
\end{equation}
consider the (relatively) open sets 
$\Sigma:=\{y\in\UU:\,u(y,0)+\eta<\psi(y)\}$,
 $\Sigma_n:=\Sigma\cap V_n$, and the set of indexes
${\mathscr I}_j:=\{\ii\in\Z^{N-1}:\,\Qijo\cap\Sigma\neq\emptyset\}$.
By the uniform continuity of $\psi$ on $V_n$
we have $\psi(y)\leq\psi(\xijo)+\eta$ for every
$y\in\cup_{\Ieps(\Sigma_n)}\Teio$ for $j$ sufficiently big.
Set $\lambdaj:=\delta_j^{1/(N-p+a)}$, define $\Teiot:=(\Teio-\xijo)/\lambdaj$
and notice that $\Teiot\subseteq B_{m-1}$ for some $m\in\N$ by (Hp 3).
Then \eqref{eq:rel-cap-unif} in Proposition~\ref{cap-prop} yields
\begin{equation}
  \label{eq:teio-unif}
\sup_{\ii\in\Z^{N-1}}|\capmu(\Teiot,B_n)- \capmu(\Teiot)|\leq\eta
\end{equation}
for all $n>m$ large enough. Let $\xiij\in\Wspo(B_{n},\mu)$ be such that
$\tilde{\xiij}\geq 1$ $\capmu$ q.e. $\Teiot$ and
$\capmu(\Teiot,B_n)=\|\nabla\xiij\|^p_{(L^p(B_n,\mu))^N}$, and let
$\zeta\in C^{\infty}_0(B_{m})$ be any function such that $\zeta\equiv 1$
on $B_{m-1}$, $\|\nabla\zeta\|_{(L^\infty(B_{m}))^N}^p\leq 2$ and $0\leq\zeta\leq 1$.

With fixed 
$n\in\N$ for which \eqref{eq:teio-unif} holds, let $(v_j)$ be 
the sequence provided by Lemma~\ref{joining} 
with $u_j\equiv u+\eta$ and $k=n$. Define 
\begin{eqnarray}
  \label{eq:recovery}
  u_j(x):=\begin{cases}
\left(1-\xiij\left(\frac{x-\xije}{\lambdaj}\right)\right)(u+\eta)_{\Cij}
+\xiij\left(\frac{x-\xije}{\lambdaj}\right)(\psi(\xijo)+\eta) &
U\cap \Bji,\,\ii\in\Ieps(\Sigma_n)\\

v_j(x) & U\setminus\cup_{{\mathscr I}_j} \Bji\\

\left(1-\zeta\left(\frac{x-\xije}{\lambdaj}\right)\right)v_j(x)
+\zeta\left(\frac{x-\xije}{\lambdaj}\right)f(x)  &
U\cap B_{m\lambdaj}(\xije),\,\ii\in{\mathscr I}_j\setminus\Ieps(\Sigma_n).
\end{cases}
\end{eqnarray}
In the definition above $\Bji$ is the set defined in \eqref{eq:bji},
and $f\in\Wsp(U,\mu)\cap L^\infty(U)$ is as in (Hp 5).

Take note that $u_j\to u+\eta$ in $L^p(U,\mu)$ since 
$\L^N(\{u_j\neq u+\eta\})\to 0$ and $U$, $u$, $f$, $\psi$ are bounded. 
Clearly, $\tilde{u}_j\geq\psi$  $\capmu$ q.e. on $\Teo$, and then
by the choice of $\xiij$, \eqref{eq:cap-scaling}, \eqref{eq:cap-sim}
and \eqref{eq:teio-unif} give
$$
\lambdaj^{-p}\int_{\R^N_+\cap\Bji}
\left|\nabla\xiij\left(\frac{x-\xije}{\lambdaj}\right)\right|^pd\mu
=\frac{\delta_j}2\capmu(\Teiot,B_{n})\leq\frac 12(\capmu(\Teio)+\delta_j\eta)
$$
for all $\ii\in\Z^{N-1}$. An analogous formula holds for the translated  
and scaled gradient of $\zeta$. 
Thus, a straightforward calculation implies
\begin{eqnarray}\label{fondam}
 \lefteqn{\hskip-1cm\FFjpsi(u_j,\om)\leq\int_U|\nabla v_j|^pd\mu+\frac 12
\sum_{\ii\in\Ieps(\Sigma_n)}\Phi(\psi(\xijo)-u_{\Cij}+\eta)
(\capmu(\Teio)+\delta_j\eta)\notag}\\&&
+2^{p-1}\|u-f\|^p_{L^\infty(U)}
\|\nabla\zeta\|^p_{L^p(B_{m})}\delta_j\#(\IIeps\setminus\Ieps(\Sigma_n))
+2^{p-1}\int_{D_j^n}|\nabla(u-f)|^pd\mu
\end{eqnarray}
where $\Phi(t)=(t\vee 0)^p$ and 
$D_j^n=\cup_{{\mathscr I}_j\setminus\Ieps(\Sigma_n)}
(U\cap B_{m\lambdaj}(\xije))$. 
By taking into account that $\L^N(D_j^n)\to 0$ as $j\to+\infty$
we may argue as in Proposition~\ref{lb} to obtain
\begin{eqnarray*}
\lefteqn{
\Gamma\hbox{-}\limsup_j\FFjpsi(u+\eta,\om)\leq
\int_U|\nabla u|^pd\mu
+\frac cn}\\&&+\frac 12\Lambda(\mathbb{E}[\gamma]+\eta)
\int_{V_n}\Phi(\psi(\xnc)-u(\xnc,0)+\eta)d\xnc
+2^{p}\|u-f\|^p_{L^\infty(U)}\Lambda\H^{N-1}(\Sigma\setminus\Sigma_n).
\end{eqnarray*}
Since $\H^{N-1}(\partial\Sigma\cup\partial\Sigma_n)=0$, by
increasing $V_n$ to $\UU$ 
we get
\begin{eqnarray}\label{fondam1}
\Gamma\hbox{-}\limsup_j\FFjpsi(u+\eta,\om)\leq
\int_U|\nabla u|^pd\mu
+\frac 12\Lambda(\mathbb{E}[\gamma]+\eta)
\int_{\UU}\Phi(\psi(\xnc)-u(\xnc,0)+\eta)d\xnc.
\end{eqnarray}
To conclude take note that $(\psi-u)\vee 0\in L^p(\UU)$
since $\FFpsi(u)<+\infty$, then there exists a positive
infinitesimal sequence $(\eta_k)$ such that
$\HH\left(\{y\in\UU:\,\psi(y)-u(y,0)=\eta_k\}\right)=0$
for every $k\in\N$. 
Moreover, $u+\eta_k\in \Wsp(U,\mu)$, being $U$ bounded,  
and it satisfies \eqref{fondam1}.
The thesis then follows by the lower semicontinuity of
$\Gamma\hbox{-}\limsup_j\FFjpsi$ as the rhs of 
\eqref{fondam1} converges to $\FFpsi(u)$ as $k\to+\infty$.


\vskip0.25cm

\emph{Step 2: $U$ unbounded.}
To remove the boundedness assumption on $U$ 
we localize the problem: for any open subset 
$A$ of $U$, $\om\in\Om$ we denote by $\FFjpsi(\cdot,\om;A)$ 
and $\FFpsi(\cdot;A)$ the functionals defined on $\Wsp(U,\mu)$ 
as $\FFjpsi$ and $\FFpsi$, respectively, with 
the domain of integration $U$ substituted with $A$.
Consider an increasing sequence $(\Ur)_{r\in\N}$ of 
bounded open Lipschitz sets in $U$ such that $\cup_r\Ur=U$, 
$B_r\cap U\subseteq\Ur\subseteq B_{r+1}\cap U$, and denote by 
$(V_n^r)_{n\in\N}$ a family of open Lipschitz subsets of $U_r$ 
such that $V_n^r\subset\subset V_{n+1}^r$ with $\cup_nV_n^r=U_r$. 
Let $\varphi_r$ be a cut-off function between $B_r$ and $B_{2r}$ with 
$\|\nabla\varphi_r\|^p_{(L^\infty(\R^N))^N}\leq 2/r^p$.  

We fix $\eta>0$ for which \eqref{eq:regularity} holds true and 
repeat for each $\Ur$ the construction of Step 1. Further, we join the 
sequence defined as in \eqref{eq:recovery} on $\Ur$ with 
the function $f$ on $\R^N\setminus\overline{U_{2r}}$. 
The sequence obtained with this construction gives the 
limsup inequality up to a vanishing error.

More precisely, with fixed $r\in\N$ and $n\in\N$ such that
\eqref{eq:teio-unif} holds, let $(u_j^r)$ be defined 
as in \eqref{eq:recovery} with $\Sigma$ and $\Sigma_n$ substituted 
by $\Sigma\cap U_{2r}$ and $\Sigma\cap V_n^{2r}$, respectively. 
Then $(u_j^r)\subset\Wsp(U_{2r},\mu)$ and 
$(u_j^r)$ converges to $u+\eta$ in $L^p(U_{2r},\mu)$. 
Define $w_j^r=\varphi_ru_j^r+(1-\varphi_r)f$,
where $f\in\Wsp(U,\mu)\cap L^\infty(U)$ is as in (Hp 5). 
Take note that $w_j^r\in\Wsp(U,\mu)$, $(w_j^r)_j$ converges to 
$(u+\eta)_r:=\varphi_r(u+\eta)+(1-\varphi_r)f$ in $L^p(U,\mu)$ 
and by definition $\tilde{w_j^r}\geq\psi$ $\capmu$ q.e. on $\UU$.
Furthermore, we have
\begin{eqnarray*}
  \lefteqn{\FFjpsi(w_j^r,\om)\leq\FFjpsi(u_j^r,\om;U_{2r})+
\FFjpsi(f,\om;U\setminus\overline{B_{2r}})}\\&&
+2^{p-1}\left(\FFjpsi(u_j^r,\om;U_{2r}\setminus \overline{B_r})
+\FFjpsi(f,\om;U_{2r}\setminus \overline{B_r})
\right)
+\frac{2^p}{r^p}\int_{U_{2r}\setminus \overline{B_r}}|u_j^r-f|^pd\mu.
\end{eqnarray*}
To estimate the rhs 
above we notice that by 
\eqref{eq:errorest} in Lemma~\ref{joining} the first and 
third terms can be dealt with as in \eqref{fondam}. 
By passing to the limsup first as $j\to+\infty$,
and then for $n\to+\infty$ we get as in \eqref{fondam1} 
\begin{eqnarray}\label{eta}
  \lefteqn{\Gamma\hbox{-}\limsup_j\FFjpsi((u+\eta)_r,\om)
\leq(1+\eta)\FFpsi(u+\eta;U_{2r})+
\int_{U\setminus\overline{B_{2r}}}|\nabla f|^pd\mu}\\&&
+2^{p-1}\left((1+\eta)\FFpsi(u+\eta;U_{2r}\setminus\overline{B_{r}})+
\int_{U_{2r}\setminus\overline{B_r}}|\nabla f|^pd\mu\right)
+\frac{2^p}{r^p}\int_{U_{2r}\setminus \overline{B_r}}|u+\eta-f|^pd\mu.
\notag
\end{eqnarray}
Arguing as in Step 1, we choose a positive infinitesimal sequence 
$(\eta_k)$ for which \eqref{eq:regularity} holds, and   
since $((u+\eta_k)_r)$ converges to $u_r$ in $L^p(U,\mu)$ as 
$k\to +\infty$, the lower semicontinuity 
of $\Gamma\hbox{-}\limsup_j\FFjpsi$ and \eqref{eta} yield
\begin{eqnarray}\label{eta2}
  \lefteqn{\Gamma\hbox{-}\limsup_j\FFjpsi(u_r,\om)
\leq\FFpsi(u;U_{2r})+
\int_{U\setminus\overline{B_{2r}}}|\nabla f|^pd\mu}\notag\\&&
+2^{p-1}\left(\FFpsi(u;U_{2r}\setminus\overline{B_{r}})+
\int_{U_{2r}\setminus\overline{B_r}}|\nabla f|^pd\mu\right)
+\frac{2^p}{r^p}\int_{U_{2r}\setminus \overline{B_r}}|u-f|^pd\mu.
\end{eqnarray}
Finally, being the rhs 
in the inequality above 
a finite measure, the lower semicontinuity of 
$\Gamma\hbox{-}\limsup_j\FFjpsi$  
gives the conclusion as $r\to+\infty$ since 
 $(u_r)$ converges to $u$ in $L^p(U,\mu)$.
\end{proof}

\begin{remark}\label{hptre}
The strong separation assumption in (Hp 3) ensures that the scaled 
obstacle sets  $(\Teio-\xijo)/\eps^{(N-1)/(N-p+a)}$ are equi-bounded and 
located in small neighbourhoods of the $x^{\ii}_j$'s.
It turns out from the proof of Theorem~\ref{main}
(see Propositions~\ref{lb} and \ref{ub}) that 
this condition can be relaxed into
\begin{itemize}

\item[{\bf (Hp 3)}$^\prime$] There exist $\e>0$ and
$\beta\in(1,(N-1)/(N-p+a)]$ such that
for all $\ii\in\Z^{N-1}$, $\om\in\Om$, and $\eps\in(0,\e)$
the sets $(\Teio-z^\ii_j(\om))/\eps^{(N-1)/(N-p+a)}$ are contained in
a fixed bounded set, for some points $z^\ii_j(\om)\in\Qij$, and
$\Teio\subseteq\xijo+\eps^\beta[-1/2,1/2)^{N-1}$.
\end{itemize}
The latter condition with $\beta>1$ is needed to ensure the validity
of the joining Lemma~\ref{joining} also in this framework.
Instead, the first condition is used when applying
Proposition~\ref{cap-prop} in the construction of the
recovery sequence in Proposition~\ref{ub}.
\end{remark}
\begin{remark}
It is clear from Remark~\ref{only-stationary}, Propositions~\ref{lb} and 
\ref{ub} that Theorem~\ref{main} still holds even dropping the ergodicity 
assumption on the $\tau_\ii$'s. The $\Gamma$-limit 
$\FFpsi:L^p(U,\mu)\times\Om\to[0,+\infty]$ being then defined as  
the functional in \eqref{eq:Glimit} with $\mathbb{E}[\gamma]$ 
replaced by the conditional expectation $\mathbb{E}[\gamma,\mathscr{I}]$.
\end{remark}
Slightly refining the argument in Step 2 above we extend the convergence 
result to the $L^p_{\mathrm{loc}}$ topology for unbounded domains.
\begin{proof}[Proof (of Theorem~\ref{main-unbdd})]
We keep using the notation introduced in Step 2 of Proposition~\ref{ub}.
The extension result in \cite[Theorem 1.1]{Chua} and Lemma~\ref{Gagliardo}
ensure that $\KpmuU$ is the domain of any $\Gamma$-cluster point.
Furthermore, the liminf inequality easily follows 
by applying Proposition~\ref{lb} to the localized functionals
$\FFjpsi(\cdot,\om;U_r)$, and then by taking the limit as $r\to+\infty$.

Instead, to get the limsup inequality for any $u\in\KpmuU$ we use
the sequence constructed in Step 2 of Proposition~\ref{ub} and 
repeat the same arguments up to \eqref{eta2}. To this aim take 
note that $\KpmuU\subset W^{1,p}_{\mathrm{loc}}(U,\mu)$.
Eventually, H\"older's inequality yields
$$
\frac{2^p}{r^p}\int_{U_{2r}\setminus \overline{B_r}}|u-f|^pd\mu\leq
2^p(\mu(B_2\setminus B_1))^{p/(N+a)}\left(
\int_{U_{2r}\setminus \overline{B_r}}|u-f|^{p^\ast}d\mu\right)^{p/p^\ast},
$$
and thus the conclusion follows as in Proposition~\ref{ub} by taking 
the limit for $r\to+\infty$ in \eqref{eta2}. 
\end{proof}

Let us now prove Corollary~\ref{bdryvalue}.
\begin{proof}[Proof (of Corollary~\ref{bdryvalue})]
The set $\Om^{\prime\prime}$ referred to in the statement 
is defined analogously to $\Om^{\prime}$. 
Hence, being $u_0+\WspoG(U,\mu)$ weakly closed, thanks to 
Proposition~\ref{lb} for all $\om\in \Om^{\prime\prime}$ we have
$$
\Gamma\hbox{-}\liminf_j(\FFepsi+\mathscr{X}_{u_0+\WspoG(U,\mu)})(u,\om)\geq
(\FF+\mathscr{X}_{u_0+\WspoG(U,\mu)})(u).
$$
Thus, given $u\in u_0+\WspoG(U,\mu)$ to conclude it suffices to verify that 
the construction of the sequence $(u_j)$ in Proposition~\ref{ub} with
$f$ there substituted by $u_0$ matches also the boundary condition
since $\overline{\Sigma}\cap\UU=\emptyset$. 
\end{proof}

\section{Generalizations}\label{genztn}

In the previous sections we have described the asymptotic behaviour of
the weighted norms on open sets $U$ subject to an obstacle condition
on part of the boundary of $U$.
In the present we discuss some generalizations of Theorem~\ref{main}.
We limit ourselves to state the results 
in these settings, since their proofs follow straightforward from the
arguments of Section~\ref{mainresult} (see also \cite{ANB}).

First, we point out that we have treated the case of the $p$-weighted
norm only for the sake of simplicity. Indeed, under only minor changes
in the proofs the same results hold for $p$-homogeneous energy densities.
Instead, the extension to non-linear energy densities having $p$-growth
seems to be more difficult. The non-linear capacitary formula introduced
by Ansini and Braides~\cite{ANB} in the deterministic setting is related
to the geometry of the scaled obstacle set. On the other hand, (Hp 1)
involves only the scaling properties of the capacity of the obstacle set,
then we are led to formulate (Hp 1)$^\prime$ below.

With fixed any $H:\R^N\to[0,+\infty)$ and $t\geq 0$ define the 
$H$-capacity of a set $E\subseteq\R^N$ by
$$
\mathrm{cap}_{H,\mu}(t,E):=\inf_{\{A \text{ open }:\,A \supseteq E\}}
\inf\left\{\int_{\R^N}H(Du)d\mu:\,u\in\Wsp(\R^N,\mu),\,
u\geq t\,\, \L^N \text{ a.e. on } A\right\}.
$$
Let $h:\R^N\to[0,+\infty)$ be a convex function such that
$$
c_1(|x|^p-1)\leq h(x)\leq c_2(|x|^p+1)
$$
for some constants $c_1$, $c_2>0$, and $h(y,x_N)=h(y,-x_N)$ for any 
$x\in\R^N$. Furthermore, put 
$H_j(x):=\eps^{\frac{(N-1)p}{N-p+a}}h\left(\eps^{-\frac{(N-1)}{N-p+a}}x\right)$
for all $x\in\R^N$.
A natural and compatible generalization of (Hp 1) is then
\begin{itemize}
\item[{\bf (Hp 1)$^\prime$.}] \emph{Capacitary Scaling:}
There exist a positive infinitesimal sequence $(\deltaj)_j$, a
function $\Phi\in C^0(\R^+)$ and
a process $\gamma:\Z^{N-1}\times\Om\to[0,+\infty)$ such that
for all $\ii\in\Z^{N-1}$ and $\om\in\Om$ we have
$$
\lim_j\mathrm{cap}_{H_j,\mu}\left(t,(\Teioe-\xijo)/\deltae^{1/(N-p+a)}\right)
=\Phi(t)\gammaio.
$$
\end{itemize}
Indeed, in case $h$ is $p$-homogeneous we have $H_j\equiv h$,
$\mathrm{cap}_{H_j,\mu}(t,E)=t^p\mathrm{cap}_{h,\mu}(1,E)$, and
thus we may take $\Phi(t)=t^p$. In the fully deterministic
setting,
i.e. $\Teio=\xijo+\deltaj^{1/(N-p+a)}T$ for some $T\subseteq\R^{N-1}$ 
for all $\om\in\Om$, $\ii\in\Z^{N-1}$ and $j\in\N$,
by assuming that $(H_j)_j$ converges pointwise to $H$ (this holds
upon extracting a subsequence by the growth conditions of $h$), we have
$\lim_j\mathrm{cap}_{H_j,\mu}(t,T)=\mathrm{cap}_{H,\mu}(t,T)$
(see \cite[Proposition 12.8]{BDF}). The continuity of
$\mathrm{cap}_{H,\mu}(\cdot,T)$ holds thanks to the local equi-Lipschitz
continuity of the $H_j$'s (which is a consequence of their 
convexity and the growth conditions of $h$).

Next we define the functional $\FF_j^h:L^p(U,\mu)\times\Om\to[0,+\infty]$ by
\begin{eqnarray}
  \label{eq:fapprox2}
  \FF_j^h(u,\om)=
  \begin{cases}
\displaystyle{\int_Uh(\nabla u)\,d\mu} &
\text{ if } u\in \Wsp(U,\mu),\,
\tildeu\geq\psi\,\,  \capmu \text{ q.e. on } \Teoe\cap\UU\\
+\infty & \text{ otherwise. }
  \end{cases}
\end{eqnarray}
The arguments by \cite{ANB} and those of Section~\ref{mainresult} then
give the following result.
\begin{theorem}\label{main2}
Assume (Hp 1)$^\prime$ and (Hp 2)-(Hp 5) hold true, $N\geq 2$, and that
$a\in(-1,+\infty)$, $p\in((1+a)\vee 1,N+a)$.

Then there exists a set $\Om^\prime\subseteq\Om$ of full probability
such that for all $\om\in\Om^\prime$ the family
$(\FF_j^h(\cdot,\om))_j$ $\Gamma$-converges in the $L^p(U,\mu)$
topology to the functional $\FF^h:L^p(U,\mu)\to[0,+\infty]$ defined by
\begin{equation}\label{eq:Glimit2}
\FF^h(u)=\int_Uh(\nabla u)\,d\mu+\frac 12
\Lambda\mathbb{E}[\gamma]
\int_{\UU}\Phi((\psi(\xnc)-u(\xnc,0))\vee 0)\,d\xnc
\end{equation}
if $u\in \Wsp(U,\mu)$, $+\infty$ otherwise.
\end{theorem}

Eventually, let us point out that similar results hold also
in case the obstacles are periodically equi-distributed inside
the open set $U$.
Clearly, conditions (Hp 1)-(Hp 5) have to be reformulated in
order to deal with the $N$-dimensional setting.
The analogue of Theorem~\ref{main2} is then an easy consequence of the
arguments of Section~\ref{mainresult} and those by \cite{ANB}.

In particular, the homogenization results in perforated open 
sets by \cite{ANB} can be extended to the ergodic setting of
Section~\ref{mainresult}, thus recovering the results of 
\cite{Caf-Mel1}, too.

\appendix

\section{ }\label{Muck} 

We show that $w(x)=|\xn|^a$ belongs to the Muckenhoupt class
$\mathcal{A}_p(\R^N)$ under a compatibility condition
between $p$ and $a$. 
\begin{lemma}\label{w-muck}
  For $p>(1+a)\vee 1$ and $a>-1$ the function $w(x)=|\xn|^a$
is in the Muckenhoupt's class $\mathcal{A}_p(\R^N)$ (see
\eqref{eq:muck} for the definition).
\end{lemma}
\begin{proof}
Let us first point out that conditions $a>-1$ and $p>1+a$
are imposed only to guarantee the local integrability of
$w$ and $w^{1/(1-p)}$, respectively.

Being $w=w(\xn)$ and even, we may restrict the supremum in
\eqref{eq:muck} to points $z=(0,z_N)$ with $z_N\geq 0$.
Define $I_\alpha(r,z_N):=\int_{B_r(z)}|\xn|^\alpha dx$
for $\alpha>-1$.
A direct integration 
yields
\begin{equation}
  \label{eq:geq2r}
  I_\alpha(r,z_N)\leq 3^{|\alpha|}\omega_Nr^N
\left(r\vee\frac{z_N}2\right)^\alpha,
\end{equation}
for $\alpha\geq 0$ and for $\alpha\in(-1,0)$ provided $z_N\geq 2r$.
Instead, in case  $\alpha\in(-1,0)$ and $z_N\in(0,2r)$ we have
$B_r(z)\subseteq B_{3r}(\underline{0})$ and again by a
direct integration we deduce
\begin{equation}
  \label{eq:leq2r}
 I_\alpha(r,z_N)\leq \frac{2N\omega_N}{1+\alpha}
\left(\int_0^1(1-t^2)^{(1+\alpha)/2}\,dt\right)(3r)^{N+\alpha}.
\end{equation}
In any case, by applying estimates \eqref{eq:geq2r} and \eqref{eq:leq2r}
above we infer
\begin{eqnarray*}
 I_a(r,z_N)\left(I_{-\frac a{p-1}}(r,z_N)\right)^{p-1}\leq
c\,r^{Np}
\end{eqnarray*}
for some positive constant $c=c(N,a,p)$. Clearly, this is equivalent to
\eqref{eq:muck}.
\end{proof}


\medskip

\end{document}